\pdfoutput=1
\documentclass[conference]{IEEEtran}
\usepackage{amsmath,amsfonts,amssymb,amsthm}

\newtheorem{theorem}{Theorem}[section]

\newtheorem{definition}[theorem]{Definition}

\newtheorem{property}[theorem]{Property}

\usepackage{multirow}
\usepackage{cite}
\usepackage{graphicx}
\usepackage{float}
\usepackage{caption}
\usepackage{subcaption}
\usepackage{bm} 
\usepackage{dsfont}
\usepackage{wasysym}
\usepackage{nameref}
\usepackage{enumerate}
\usepackage{nomencl}
\usepackage{etoolbox}
\usepackage[symbol]{footmisc}

\begin{document}
	\bstctlcite{MyBSTcontrol}
	\title{Chance-Constrained Ancillary Service Specification \\for Heterogeneous Storage Devices}
	
	\author{\IEEEauthorblockN{Michael P. Evans, David Angeli$^*$ and Goran Strbac}
		\IEEEauthorblockA{Department of Electrical and Electronic Engineering\\
			Imperial College London, UK\\
			\{m.evans16, d.angeli, g.strbac\}@imperial.ac.uk}
		\and
		\IEEEauthorblockN{Simon H. Tindemans}
		\IEEEauthorblockA{Department of Electrical Sustainable Energy\\
			TU Delft, Netherlands \\
			s.h.tindemans@tudelft.nl}
	}
	
	\maketitle
	\footnotetext[1]{Dr David Angeli is also with the Department of Information Engineering, University of Florence, Italy}
	
	\begin{abstract}
		We present a method to find the maximum magnitude of any supply-shortfall service that an aggregator of energy storage devices is able to sell to a grid operator. This is first demonstrated in deterministic settings, then applied to scenarios in which device availabilities are stochastic. In this case we implement chance constraints on the inability to deliver as promised. We show a significant computational improvement in using our method in place of straightforward scenario simulation. As an extension, we present an approximation to this method which allows the determined fleet capability to be applied to any chosen service, rather than having to re-solve the chance-constrained optimisation each time.
	\end{abstract}
    
    \vspace{2mm}
	
	\section{Introduction}
	Recent years have seen a significant increase in the proliferation of energy storage devices onto electricity networks. This is in large part due to their potential to replace conventional generation as electricity grids are decarbonised \cite{WEC}. As part of this transition, increasing numbers of system operators offer tenders for the provision of ancillary services; in this paper we focus on balancing services, although our analysis would be applicable to other active power services as well. The majority of these services are specified in terms of time durations over which a ramp in power or a constant power output should be delivered. When considered for provision by storage, they are therefore often taken to be shaped delivery requirements, as in \cite{Greenwood2017}. We place our analysis within such a framework, and in particular consider request profiles in the shape of pulses or trapezoids. We point out here, however, that our analysis is general and not limited to any particular shape.
	
	Many prior studies have considered the application of storage to the provision of ancillary services. A particularly promising field is that of vehicle-to-grid (V2G), in which electric vehicles (EVs) are able to provide system support when connected to the network \cite{Ota2012,Sortomme2012,Han2011,Kisacikoglu2010}. This has been investigated for a range of services, including spinning reserve \cite{Ota2012,Sortomme2012}, regulation \cite{Sortomme2012,Han2011} and reactive power compensation \cite{Kisacikoglu2010}. Crucially, the ability for an EV to participate is subject to its availability, which tends to be assumed stochastic. We incorporate this form of stochasticity into our analysis. Moreover, in these settings a common requirement is that the provision of an ancillary service must not reduce the quality of the primary service (e.g. transportation) below some user-defined threshold. A similar framework might apply to uninterruptible power supplies (UPSs) that are kept as a backup in case of power outage. Other device configurations instead include system service provision within their core tasks, as is common among home energy management systems including storage \cite{Zhou2016}. Such a system might, for example, be designed so as to trade off supporting the grid against maximising self-consumption of electricity generated from rooftop solar panels \cite{Moshovel2015}. Wu \textit{et al.} \cite{Wu2018} combine the strands of V2G and home energy management by considering a system in which the only storage device is the household EV.
	
	It is clear from these applications that a means to determine \textit{a priori} the maximum service provision capability of a fleet of heterogeneous storage devices would be a valuable asset to an aggregator. This would allow them to determine the magnitude of a given shape of delivery requirement that they could offer to the system operator. When considering a single storage device, it is trivial to deduce its maximum ancillary service capability, simply by considering the power and energy limits of the device. The same cannot be said of a heterogeneous fleet of devices, however. In our prior work \cite{Evans2018}, we presented a transform through which one can represent the aggregate flexibility of a device pool, termed the $E$-$p$ transform. An alternative aggregate fleet representation was developed independently by Cruise and Zachary \cite{Cruise2018}. For an aggregate representation of this form to exist requires that the problem be restricted to one of discharge only, as was discussed in detail in \cite{Evans2018}.
	
	This paper explores in greater depth how the $E$-$p$ transform can be used to specify ancillary service capability. In particular, we move beyond deterministic settings. Chance constraints are commonly applied to the dispatch of storage under uncertainty \cite{Akhavan-Hejazi2018,Zhao2018,Xu2018}, and we utilise these to ensure that the probability of unsuccessful delivery of the procured service lies below a user-defined threshold. This involves a Monte Carlo simulation in which the $E$-$p$ transform of \cite{Evans2018} is used as a binary feasibility check on each drawn trace. We point out here that an alternative method for such a check would be a procedural simulation in time spanning the total duration of the delivery requirement. We therefore investigate the speed-up achieved in using the $E$-$p$ transform for such a check. In prior work we presented a policy, in both continuous \cite{Evans2017} and discrete \cite{Evans2018a} time, with a one-to-one mapping between feasible requests and those that will be met under the implementation of the policy. It will therefore be this that we will use to perform a procedural check.
	
	\section{Mathematical Formulation}
	\label{sec: math}
	\vspace{-1mm}
	\subsection{Problem description}	
	We denote the constituent units of the heterogeneous fleet as $D_i,\ i=1,...,n$, each with extractable energy $e_i(t)$ subject to the assumed physical constraint $e_i(t)\geq0$. We choose as our control input $u_i(t)$ the power extracted from each device, measured externally so as to take into account any inefficiencies present, leading to the integrator dynamics $\dot{e}_i(t)=-u_i(t)$. For the reasons mentioned in the Introduction, we restrict ourselves to discharging operation, with device-specific limit $\bar{p}_i$, so that $u_i(t)\in[0,\bar{p}_i]$. We choose as our state variable the \textit{time-to-go}, defined for each device as $x_i(t)\doteq e_i(t)/\bar{p}_i$. We then stack state, input and maximum power values across devices as follows:
	\begin{alignat}{3}
	x(t)&\doteq \begin{bmatrix}
	x_1(t)\ \dots\ x_n(t)
	\end{bmatrix}^T,\\
	u(t) &\doteq \begin{bmatrix}
	u_1(t)\ \dots\ u_n(t)
	\end{bmatrix}^T,\\
	\bar{p}&\doteq\begin{bmatrix}
	\bar{p}_1\ \dots\	\bar{p}_n
	\end{bmatrix}^T,
	\end{alignat}
	allowing us to rewrite the dynamics in matrix form as $\dot{x}(t)=-P^{-1}u(t)$, in which $P\doteq\text{diag}(\bar{p})$. We also form the product set of the power constraints, $\mathcal{U}_{\bar{p}}\doteq[0,\bar{p}_1]\times[0,\bar{p}_2]\times...\times[0,\bar{p}_n]$,	so that our constraints can be compactly written as $u(t)\in\mathcal{U}_{\bar{p}}$ and $x(t)\geq0$. We will be interested in the ability of the fleet to meet a request profile, which we denote by $P^r\colon[0,+\infty)\mapsto[0,+\infty)$. We say that any profile which the fleet can satisfy without violating any constraints is \textit{feasible}, and we define the set of such signals as follows:
	\begin{definition}
		\label{def:feasible}
		The set of feasible request profiles, for a system with maximum power vector $\bar{p}$ and initial state $x=x(0)$, is defined as
		\begin{equation*}
		\begin{aligned}
		\mathcal{F}_{\bar{p},x} \doteq \big\{&P^r(\cdot) \colon \exists u(\cdot),\ z(\cdot)\colon \forall t\geq 0,\ 1^Tu(t)=P^r(t),\\
		&u(t)\in \mathcal{U}_{\bar{p}},\ \dot{z}(t)=-P^{-1}u(t),\ z(0)=x,\ z(t) \geq 0 \big\}.
		\end{aligned}
		\end{equation*}
	\end{definition}
	\noindent This then allows us to transform a determination of request feasibility into an evaluation of membership of this set. 
	
	\subsection{$E$-$p$ transform}
	The main tool that we will use in our analysis is the $E$-$p$ transform of \cite{Evans2018}. We therefore reproduce the following definitions here:
	\begin{definition}
		Given a request profile $P^r\colon [0,+\infty)\to[0,+\infty)$, we define its E-p transform as the following function:
		\begin{equation}
		E_{P^r}(p)\doteq\int_0^\infty \textnormal{max}\big\{P^r(t)-p,0\big\}dt,
		\end{equation}
		interpretable as the energy required above any given power rating, $p$.
	\end{definition}
	\begin{definition}
		\label{definition: capacity}
		We define the \textbf{capacity} of a system to be the $E$-$p$ transform of the worst-case request profile that it can meet, i.e. $\Omega_{\bar{p},x}(p)\doteq E_{R}(p)$, where $R(\cdot)$ is defined as
		\begin{equation}
		R(t)\doteq\sum_{i=1}^n\bar{p}_i[H(t)-H(t-x_i)],
		\end{equation}
		in which $H(\cdot)$ denotes the Heaviside step function.
	\end{definition}
	
	\begin{property}
		\label{the:E_p}
		A request profile $P^r(\cdot)$ is feasible if, and only if, its $E$-$p$ transform is dominated by the capacity curve of the system, i.e. $E_{P^r}(p)\leq \Omega_{\bar{p},x}(p)\ \ \forall p\iff P^r(\cdot)\in\mathcal{F}_{\bar{p},x}$.
	\end{property}
	
	\noindent This property can be interpreted to mean that the capacity curve represents the full set of request profiles that a fleet can satisfy. As discussed in detail in \cite{Evans2018}, one can therefore use this curve to determine the maximum ancillary service capability of the fleet. 
	
	In this paper, we use the following bisection method to find the maximum feasible service magnitude for a given capacity curve. Start with the interval $[0,\ \sum_{i=1}^n\bar{p}_i]$, which is known to contain the maximum feasible magnitude. Then, continually bisect and update this interval until its width is less than a chosen tolerance. Finally, as we are interested in the largest feasible magnitude, take as the solution the lower bound. 
	
	\section{Application to Stochastic Settings}
	\label{sec:stoch}
	\subsection{Problem of interest}
	We take the point of view of an aggregator that wants to sell grid support services to a system operator. More specifically, we restrict these to balancing services under supply-shortfall conditions, so that our delivery requirements are strictly positive. We consider the setting in which devices have stochastic availability, which results in a random initial state vector and therefore a random capacity curve. We investigate the scenario in which an aggregator offers a service to a grid operator of a predefined shape and duration, but where they are free to specify the magnitude, $m$, of that service; and they are incentivised to offer as large a value as they can. We therefore aim to find the maximum such magnitude that the aggregator can sell, subject to a chance constraint on feasibility. We model their risk aversion via a level $c$, encoding the maximum probability of being unable to supply the contracted service that they are willing to accept. Their task, for a given $P^r(\cdot)$, is then to solve the following chance-constrained optimisation problem:
	\begin{equation}
	    \begin{aligned}
	    \max_m\ &m P^r(\cdot) \\
	    \textnormal{s.t. } \textnormal{Pr}\ [&m P^r(\cdot)\in\mathcal{F}_{\bar{p},x}]\geq 1-c,
	    \end{aligned}
	    \tag{$\star$}
	    \label{eq:optimisation problem}
	\end{equation}
	in which Pr$[\cdot]$ denotes the probability operator. We consider scenarios in which a probability distribution of device availabilities is known. Despite this, however, we do not have an analytical expression for the probability distribution of interest. Instead, we are able to sample from this distribution, and we implement Monte Carlo sampling to find an estimate of the solution to \eqref{eq:optimisation problem}.
	
	There are many scenarios to which this might apply in reality. For example, availability of EVs to participate in V2G services is clearly limited by whether or not they happen to be plugged in when the request is received. Similarly, UPSs might be allowed to offer regulation services, but this would have to be curtailed in the event that they were called upon to deliver backup power. Under any such scenario with stochastic device availability, a chance-constrained capability assessment would tell the aggregator what magnitude of service provision to offer the system operator. We apply the $E$-$p$ transform to this problem for the following reasons. Not only does it offer a speed-up over procedural simulation, as we will demonstrate, but it also offers intuition and insights into the factors affecting fleet capability. This enables us to find routes to computationally efficient approximations.
	
	\vspace{-1.5mm}
	\subsection{Monte Carlo procedure}
\begin{figure}[h]
		\centering
		\vspace{-8mm}
		\centerline{\includegraphics[width=0.9\columnwidth]{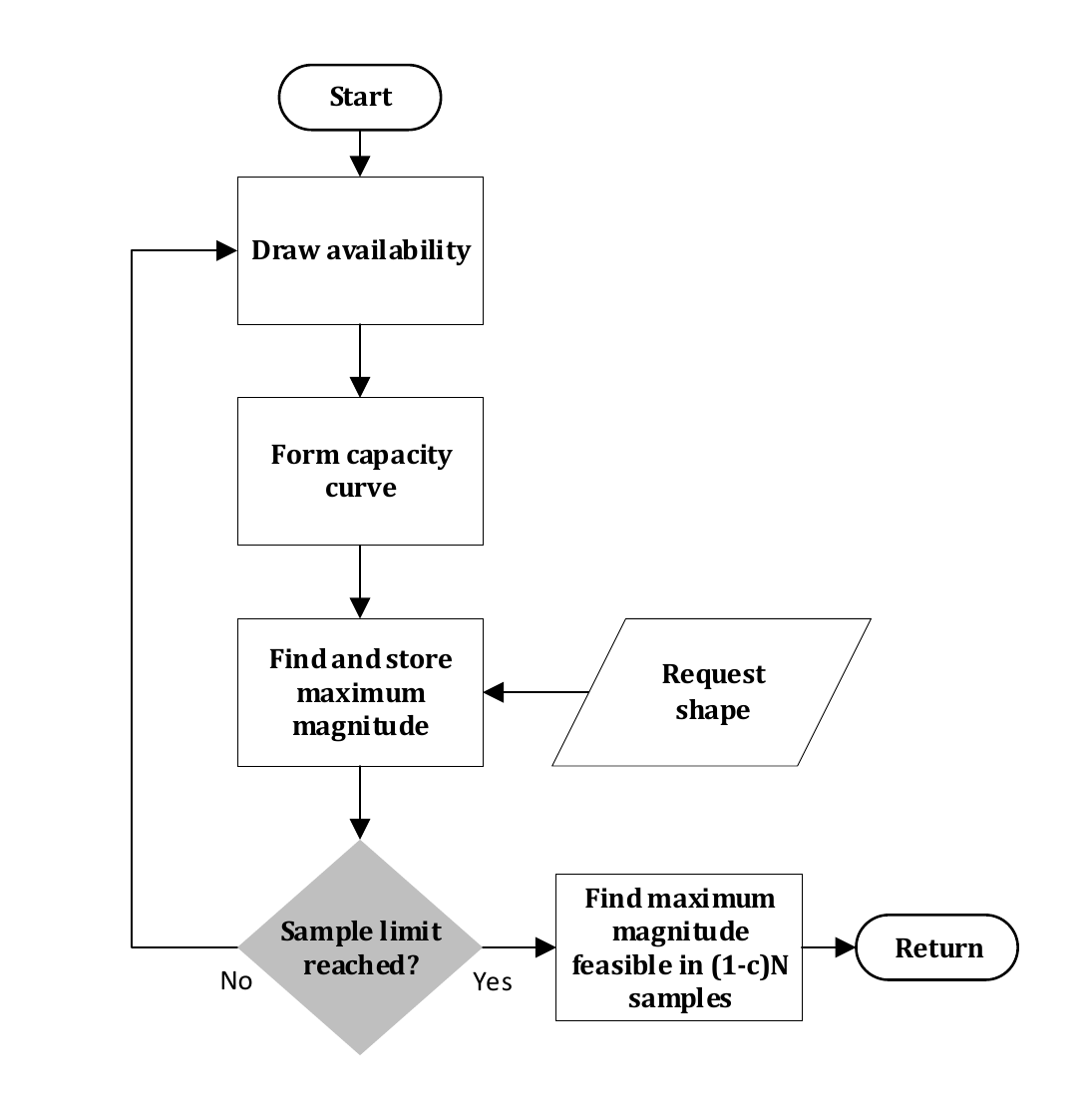}}
		\vspace{-4mm}
		\caption{The Monte Carlo procedure}
		\label{fig:flowchart}
	\end{figure}

	As mentioned previously, for a given maximum power and state vector, we can find the maximum service capability of the fleet using its capacity curve. In our problem framework, the maximum power vector is predefined, but the state vector is allowed to vary based on device availability. We denote by $a_i\in\{0,1\}$ the availability of device $D_i$ and stack availabilities across devices as $a\doteq\begin{bmatrix}a_1,...,a_n\end{bmatrix}^T$. The sample-state corresponding to a given realised availability is then $x^s=a\ {\circ}\ x$, the Hadamard product of the availability and the (full-availability) state. For a given full-availability state, our Monte Carlo procedure utilises this as follows. For each sample, we draw an availability vector $a$ and compute $\Omega_{\bar{p},x^s}$. Using this, we then find the maximum feasible magnitude of the chosen service. We repeat this $N$ times, then take the maximum magnitude which would have resulted in feasibility in $(1-c)N$ of the samples. This procedure can be seen diagrammatically in Figure~\ref{fig:flowchart}.
	
	\vspace{-5mm}
	\subsection{Quantile approximation}
\begin{figure}[h]
		\centering
		\vspace{-8mm}
		\centerline{\includegraphics[width=0.9\columnwidth]{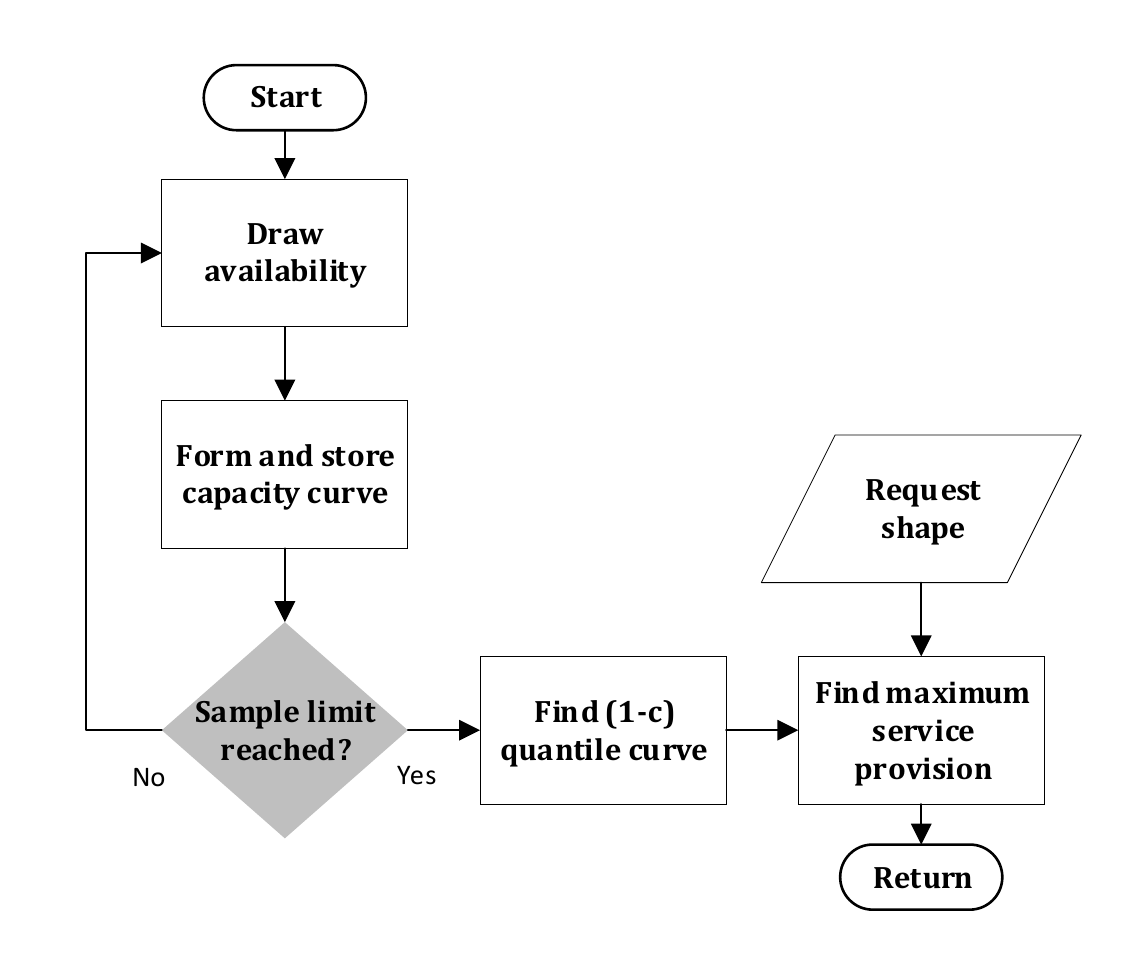}}
		\vspace{-3.5mm}
		\caption{The quantile approximation procedure}
		\label{fig:flowchart2}
	\end{figure}

	We here present an approximation to the above procedure, as follows. We draw $N$ samples and each time form the corresponding capacity curve. We then determine the single curve such that $(1-c)N$ of the observed capacity curves lie below it, point-wise in $p$. The single curve formed is then used as a deterministic capacity. This procedure can be seen in Figure~\ref{fig:flowchart2}. We point out here that this is an optimistic approximation, because pointwise dominance of $(1-c)N$ capacity curves is a less stringent condition than dominance over the whole domain of $p$. If this approximation is observed to be reasonably accurate, however, it might still be worthwhile. By producing a single curve, it would enable the determination of the maximum magnitude of \textit{any} ancillary service, without needing to compose new samples each time. In addition, optimisation with respect to a quantile curve should be much simpler than optimisation with respect to an ensemble.
	
	\vspace{-3mm}
	\section{Numerical Results}
	\label{sec:num}
	\vspace{-1mm}
	\subsection{Scenario studied}
	A case study illustrates the procedures described above as follows. We choose the fleet size to be 500 devices. We draw maximum power and energy ratings for each device from device-independent log-normal distributions, based on the assumption that relative variation in these values across the fleet is more important than absolute variation. Inspired by typical EV parameters, we assign 450 (90\%) of the devices to have maximum power values with a mean of 3.3 kW and a standard deviation of 1 kW, each of the log-normal distribution, and the remaining 50 devices to have a rating of 50 kW. This represents the partial access to rapid chargers, under the assumption that all chargers have V2G capability equal to their rated power. We then choose the current energy values of the (entire) fleet to have a mean of 40 kWh and a standard deviation of 10 kWh, again of the log-normal distribution.
	
	\subsection{Ancillary service specification procedure}
\begin{figure}[h]
		\centering
		\begin{subfigure}[l]{\columnwidth}
			\centering
			\includegraphics[width=0.61\columnwidth]{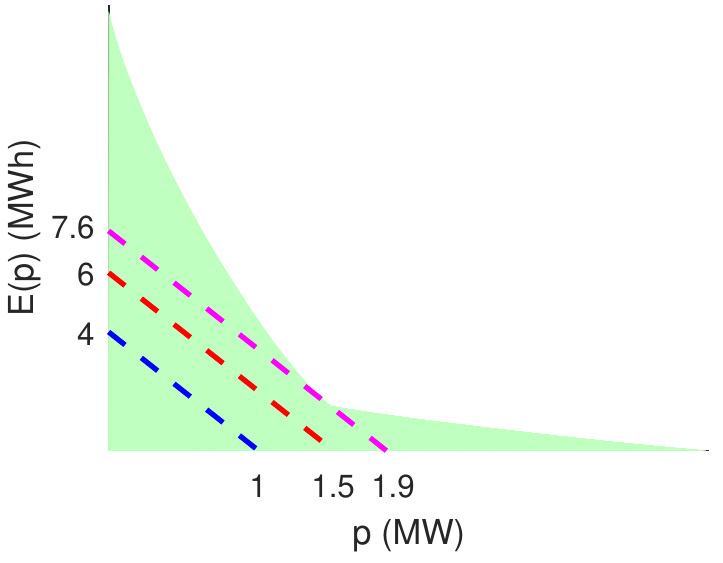}
			\caption{Use of the $E$-$p$ transform to determine the maximum-magnitude 4 h pulse that the fleet can provide. The feasible region is shaded in green and the requests are represented by dashed lines.}
			\label{fig:pulse_Ep}
		\end{subfigure}%
		\quad
		\begin{subfigure}[l]{\columnwidth}
			\centering
			\vspace{5mm}
			\includegraphics[width=0.55\columnwidth]{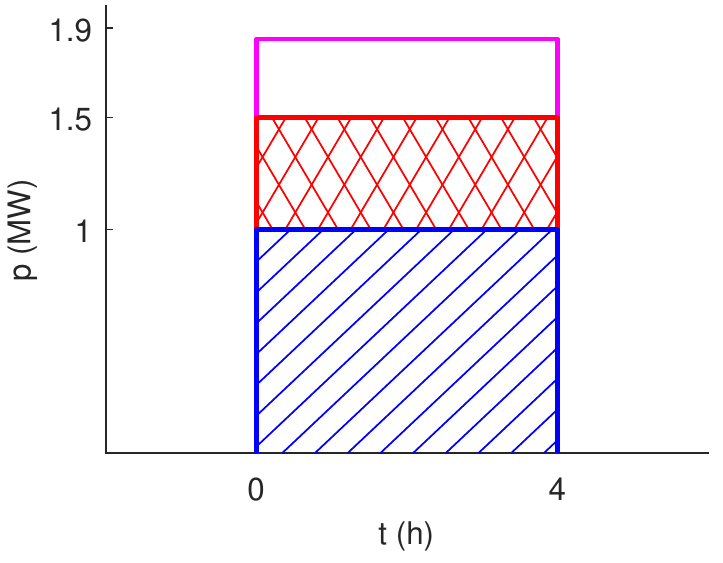}
			\caption{The three 4 h pulses considered. To enable a direct comparison, the unfilled magenta pulse lies behind the cross-hatched red pulse, which in turn lies behind the hatched blue pulse.}
			\label{fig:pulse_signal}
		\end{subfigure}
		\caption{Use of the $E$-$p$ transform to determine the largest feasible pulse magnitude of 4 h duration. The colour of each $E$-$p$ curve in Figure~\ref{fig:pulse_Ep} matches that of the corresponding pulse signal in Figure~\ref{fig:pulse_signal}.}
	\end{figure}

	Initially, we consider the capability of the fleet when all devices are fully available. In this case all parameter values are known explicitly and therefore the problem is deterministic. For this initial example we consider a service to deliver a simple pulse with a duration of 4 h, and find the maximum magnitude pulse of this form that the fleet can provide. The $E$-$p$ curve corresponding to a signal of this type will be linear, with gradient equal to the negative duration and $p$-intercept equal to the pulse magnitude. Therefore this task involves finding the straight line of gradient -4, with the maximum $p$-intercept out of those that are contained within the feasible region. This procedure can be seen in Figure~\ref{fig:pulse_Ep}, with the corresponding pulses shown in Figure~\ref{fig:pulse_signal}.
	
	\vspace{-3mm}
	\subsection{Trapezoidal delivery requirements}
\begin{figure}[h]
		\centering
		\begin{subfigure}[l]{\columnwidth}
			\centering
			\includegraphics[width=0.675\columnwidth]{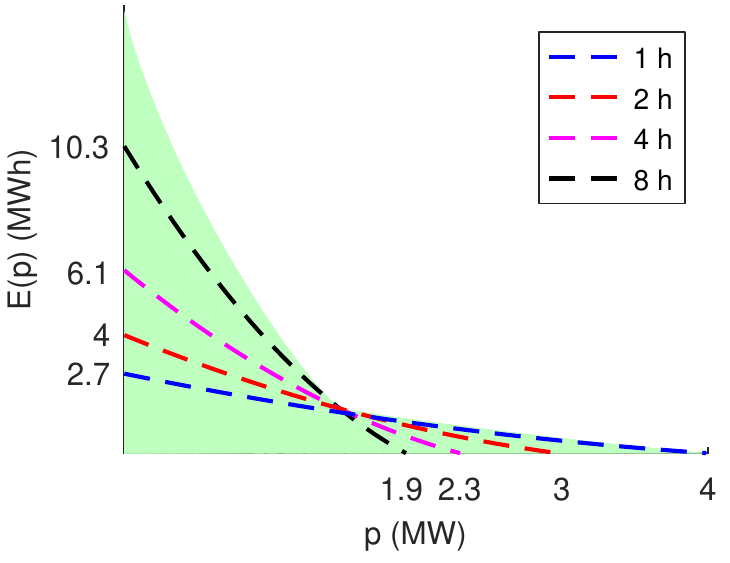}
			\caption{Use of the $E$-$p$ transform to determine the maximum-magnitude trapezoid that the fleet can provide, across a range of total durations. The feasible region is shaded in green and the requests are represented by dashed lines, with their total durations given in the figure legend.}
			\label{fig:trapezoid_Ep}
		\end{subfigure}%
		\quad
		\begin{subfigure}[l]{\columnwidth}
			\centering
			\vspace{5mm}
			\includegraphics[width=0.65\columnwidth]{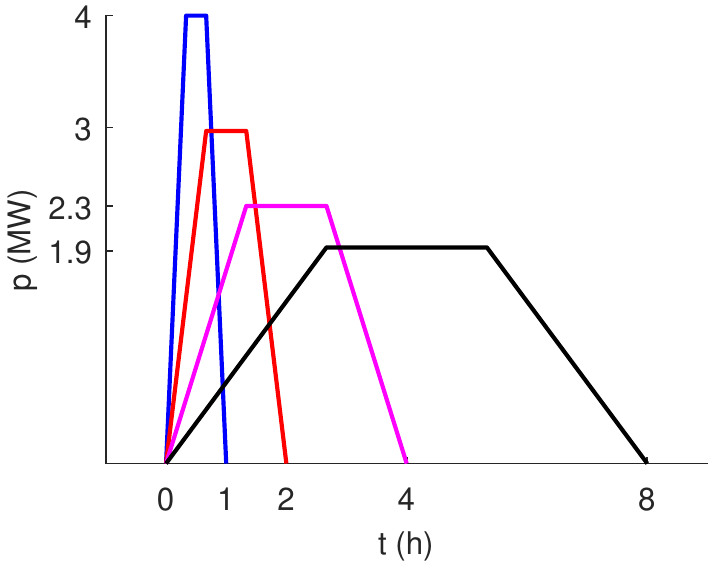}
			\caption{The trapezoidal requests considered.}
			\label{fig:trapezoid_signal}
		\end{subfigure}
		\caption{Use of the $E$-$p$ transform to determine the largest-magnitude feasible trapezoid, across a range of total durations. The colour of each $E$-$p$ curve in Figure~\ref{fig:trapezoid_Ep} matches that of the corresponding pulse signal in Figure~\ref{fig:trapezoid_signal}.}
	\end{figure}

	Having demonstrated how one can use the $E$-$p$ transform to specify maximum ancillary service provision, we now perform this task on a less straightforward service. We consider a trapezoidal signal in which the total duration is split equally among ramping up at a fixed rate, maintaining constant power and ramping down to zero at the (negative) same fixed rate. We find the largest magnitude of such a signal that the fleet can provide, across a range of total durations. Note that, as in the previous example, the magnitude of this delivery requirement is given by the $p$-intercept of each curve. The results of this can be seen in Figure~\ref{fig:trapezoid_Ep}, with the corresponding request profiles shown in Figure~\ref{fig:trapezoid_signal}.
	
	\subsection{Chance-constrained optimisation}
	We now consider the case in which device availability is stochastic, varying according to a known probability distribution. We assume an independent Bernoulli distribution for each device, with a 60\% probability that the device is available. We implement a Monte Carlo simulation to determine the maximum magnitude 2 h trapezoidal signal that the fleet can now provide, across chance constraint probability levels, $c$, of 50\%, 10\% and 1\%. We apply the procedure of Figure~\ref{fig:flowchart}, with the number of samples, $N$, set to $10^4$.	The results of this experiment are shown in Table~\ref{tab: M values} and plotted as dashed $E$-$p$ curves in Figure~\ref{fig:quantile_approx}. 95\% bootstrap confidence intervals were observed to be smaller than the significant digits shown.

	\vspace{-2mm}
	\subsection{Quantile approximation}
\begin{figure}[h]
		\centering
		\vspace{-3mm}
		\includegraphics[width=0.78\columnwidth]{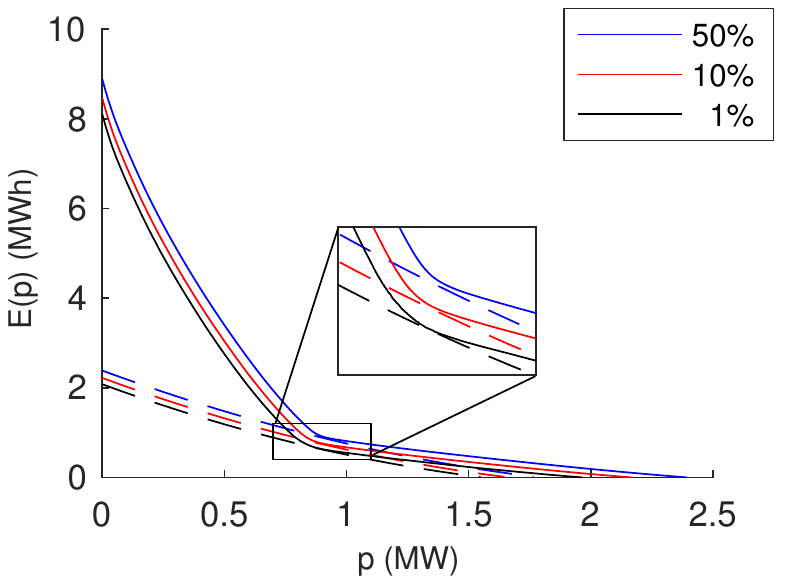}
		\caption{Approximating the above results using a single quantile curve, across a range of probability levels, $c$ (shown in the figure legend). In each case, the accurate maximum $E$-$p$ curve is represented by a dashed line, with the approximate capacity curve represented by a solid line of the same colour.}
		\label{fig:quantile_approx}
	\end{figure}

	We now investigate the approximation to the chance-constrained result described in Section~\ref{sec:stoch}. We use the same fleet as in the previous experiment, including the same distribution on availabilities of devices, and we again set $N=10^4$.	The results of this can be seen in Table~\ref{tab: M values}, and the single-curve approximations are plotted as solid lines in Figure~\ref{fig:quantile_approx}. Again, 95\% bootstrap confidence intervals were observed to be smaller than the significant digits shown. Interestingly, this experiment returns a magnitude of provision that differs from the accurate value by less than 1\%, suggesting that this quantile technique may be a worthwhile approximation to make. This is corroborated by Figure~\ref{fig:quantile_approx}, where for each probability level there is a very good match (i.e. a small gap) between the accurate request $E$-$p$ curves and the approximate capacity curves. We reiterate, however, that this estimate is not robust; we know from the preceding results that each approximated magnitude would lead to a greater probability of failure than the aggregator is willing to accept.
	
	\renewcommand{\arraystretch}{1.2}
	\begin{table}[h]
		\centering
		\begin{tabular}{|c|c|c|c|}
			\hline
			\multirow{2}{*}{$c$ value} & Accurate  & Approximated & Relative\\
			 & magnitude (MW) & magnitude (MW) & Error\\
			\hline
			50\% & 1.79 & 1.80 & 0.98\% \\
			10\% & 1.66 & 1.67 & 0.35\% \\
			\hspace{1.5mm}1\% & 1.55 & 1.56 & 0.75\% \\
			\hline
		\end{tabular}
		\caption{The magnitude values determined by each procedure, and the approximation errors.}
		\label{tab: M values}
	\end{table}
	\renewcommand{\arraystretch}{1}
	
	\vspace{-3mm}
    
    \vspace{2mm}
    \subsection{Timing comparison}
	We now investigate the speed-up achieved in using the $E$-$p$ transform to determine the maximum service magnitude for a given sample-state, in place of straightforward time-sequential simulation using the discrete time policy of \cite{Evans2018a}. We assume that the trapezoidal request is implemented as a piecewise constant signal of period 1 min, and for increased efficiency configure the routine so that it terminates as soon as the current value of the request cannot be met. A binary feasibility check is called once per iteration of the bisection procedure, but the number of iterations required depends on the sampled state. For a direct comparison, therefore, we find the average time taken across the bisection procedure, over the same $10^4$ sampled states in each case. For the $E$-$p$ case we include formation of the capacity curve within this, since it is not required during the straightforward simulation. The results of this comparison can be seen in Table~\ref{tab: speed-up}, where it can be seen that the $E$-$p$ transform enables a 2.6x improvement over straightforward simulation. This experiment was performed in MATLAB on a workstation with an Intel Xeon processor and 64 GB of memory. 
	
	\renewcommand{\arraystretch}{1.2}
	\begin{table}[h]
		\centering
		\begin{tabular}{|l|c|}
			\hline
			Method & Average run-time (ms)\\
			\hline
			Discrete time optimal policy & 3.42 \\
			$E$-$p$ transform & 1.33  \\
			\hline
		\end{tabular}
		\caption{The total computation time using each method.}
		\label{tab: speed-up}
	\end{table}
	\renewcommand{\arraystretch}{1}

    \vspace{-3mm}
	\section{Conclusions and Future Work}
	\label{sec:conc}
	This paper has demonstrated use of the $E$-$p$ transform to determine request feasibility, in both deterministic and stochastic settings. We have shown how this approach can be used to specify maximum service provision for a fleet of heterogeneous storage devices, and demonstrated that the speed-up achieved as compared to procedural simulation can be significant. In deterministic settings the maximum service capability takes the form of a single curve, whereas in stochastic settings chance constraints can be utilised to ensure that the probability of successful delivery lies above a chosen limit. We have also investigated an approximation to this chance-constrained analysis. By forming a single capacity curve, this enables chanced-constrained optimisation to be performed across any chosen service.
	
	In future work the authors plan to develop an analytical basis for this approximation, investigate its accuracy when supplied with fewer samples and further explore its applicability. They also intend to embed these procedures into larger optimisation routines.
	
	\bibliographystyle{IEEEtran}
	\bibliography{bib_new2,biblio_manual3}
\end{document}